\documentclass[12pt,a4paper,twoside,reqno]{amsart}

\newtheorem{corollary}{\sc Corollary}[section]
\newtheorem{theorem}{\sc Theorem}[section]
\newtheorem{lemma}{\sc Lemma}[section]
\newtheorem{proposition}{\sc Proposition}[section]
\newtheorem{remark}{\sc Remark}

\usepackage{amsmath}
\usepackage{amsfonts}
\usepackage{amssymb}
\usepackage{amstext}
\usepackage{amsbsy}
\usepackage{epsf}

\let\hip=\Hip

\newcommand{\R}{\mathbb R}
\newcommand{\M}{\mathbb M}

\newcommand{\wt}{\widetilde}

\let\h=\hip

\let\re=\R

\def\s{\mathbb  S}

\def\leqs{\leqslant}

\def\rmd{\mathop{\rm d\kern -1pt}\nolimits}
\def\rme{\mathop{\rm e\kern -1pt}\nolimits}

\DeclareMathOperator{\dist}{dist}

\DeclareMathOperator{\ext}{ext}

\def\bel{ \medskip
 \centerline{$ \ast \hbox to 1.0cm{}\ast \hbox to 1.0cm{}\ast $}
}

\def\longerrightarrow{-\kern-5pt\longrightarrow}

\def\star{\lower 1pt\hbox{*}}
\def \nulset {
\raise 1pt\hbox{ \hskip -3pt$\not$\kern -0.2pt \raise
.7pt\hbox{${\scriptstyle\bigcirc}$}}}

\newcommand{\hd}{\mathbb{H}^2}
\newcommand{\sd}{\mathbb{S}^2}
\newcommand{\hi}[1]{\mathbb{H}^#1}
\newcommand{\m}[1]{\mathbb{R}^#1}
\newcommand{\ch}{\cosh}
\newcommand{\sh}{\sinh}

\newcommand{\pain}{\partial_{\infty}}
\newcommand{\var}{\varepsilon}
\newcommand{\ov}[1]{\overline{#1}}

\DeclareMathOperator{\argsh}{arcsinh}
\DeclareMathOperator{\argch}{arcosh}

\def\lra{\longrightarrow}

\begin{document}

\title[existence and uniqueness]
{Uniqueness of $H$-surfaces in $\hd \times \R$, $\vert H\vert \leq 1/2$,
with boundary one or two
parallel horizontal circles}

\author[ B. Nelli, R. Sa Earp,\\ W. Santos and E. Toubiana]
{\scshape B. Nelli, R. Sa Earp,\\ W. Santos and E. Toubiana}

\address{Dipartimento di Matematica,
Universit\'a di L'Aquila, \newline
via Vetoio - Loc. Coppito,
67010 (L'Aquila), Italy}\email{nelli@univaq.it}
 \address{Departamento de Matem\'atica,
  Pontif\'\i cia Universidade Cat\'olica do Rio de Janeiro, Rio de Janeiro,  22453-900 RJ,
 Brazil }\email{earp@mat.puc-rio.br}
\address{Universidade Federal do Rio de Janeiro, Instituto de
Matemática, Av. Brigadeiro Trompowsky, s/nº Cidade Universitária,
Ilha do Fund\~ao, Caixa Postal 68530, 21945-970 Rio de Janeiro RJ,
Brazil} \email{walcy@im.ufrj.br}
\address{Institut de Math{\'e}matiques de Jussieu, Universit{\'e} Paris
VII, Denis Diderot, Case 7012,
         2 Place Jussieu,
         75251 Paris Cedex 05, France}
\email{toubiana@math.jussieu.fr}
\thanks{The second and the third authors would like to thank CNPq,  PRONEX of
Brazil and Accord Brasil-France, for
partial financial support}



\maketitle

\begin{abstract}
 We prove that a $H$-surface $M$ in $\hd \times \R$,
$H\leq\frac{1}{2}$, inherits the symmetries of its boundary $\partial M,$ when
$\partial M$ is either a horizontal curve with curvature greater than one or
two  parallel horizontal curves  with curvature greater than one, whose
distance is greater or equal to $\pi.$
Furthermore we prove that the asymptotic boundary of a surface with  mean
curvature  bounded away from zero consists of parts of straight lines, provided
 it is sufficiently regular.
\end{abstract}

\section{Introduction }

An old question in classical Differential Geometry in Euclidean space concerns
the influence of the boundary on the behavior of a $H$-surface.
A similar question can also be asked when the ambient space is a homogeneous
3-manifold.

 In this paper we consider this problem when the ambient space is the product
$\hd \times \R$ and the boundary of the $H$-surface $M$ is  either a horizontal
curve (Theorem
\ref{oneboundary}) or two parallel  horizontal curves (Theorem
\ref{smalldistancetheorem})
with curvature greater than one. By parallel curves we mean
congruent up to a vertical translation.
The main result is that $M$ inherits the
symmetries of its boundary.
In particular, if the boundary curve
is a horizontal circle
 or two parallel horizontal circles with distance greater or equal to $\pi$,
then $M$ is rotational.

We recall the principal related results in Euclidean and
hyperbolic 3-space.

  It has been conjectured that a connected compact embedded $H$-surface in
$\m 3$  with boundary a round circle is spherical. Of course, this conjecture
can be posed in hyperbolic space $\hi 3$ and in the product space
$\hd \times \R$. In $\m 3$ this is still an open problem. F. Braga Brito, W.
Meeks, H. Rosenberg and R. Sa Earp proved the conjecture provided
the surface is transverse to the plane containing the
circle, along the circle (\cite{BMRS}). F. Braga Brito and R. Sa Earp
proved that, if the radius of the circle and the mean curvature are
equal to one, then the conjecture is true assuming only that the surface is
immersed (\cite{FSE1}). In fact they also proved analogous
characterizations of a spherical cap for $f$-surfaces (special Weingarten
surfaces) of disk type (\cite{FSE2}). We remark that N. Kapouleas
announced examples of immersed $H$-surfaces with genus $g\geq 3$, with boundary
a
circle (\cite{Ka}).

 In hyperbolic space L. Barbosa and
R. Sa Earp proved the following sharp result: if a compact  connected immersed
surface in $\hi 3$, with
boundary a round circle, has constant mean curvature smaller or equal to one,
then it is totally umbilical (\cite{BSE1}, \cite{BSE2} and \cite{BS}). B. Nelli
and
H. Rosenberg proved the same result in the embedded case (\cite{Ne-Ro}). In
\cite{SET2}, R. Sa Earp and E. Toubiana generalized the above result of
\cite{Ne-Ro} to $f$-surfaces in
hyperbolic 3-space satisfying $f^2\leq 1$.
They also
proved the following: if $M$ has
constant mean curvature one, is embedded into $\hi3$ with
$\partial M=C_1 \cup C_2$,
where $C_1$ and $C_2$  are two parallel circles and
$\dist (C_1,C_2)$ great enough, then $M$ is a piece of a catenoid
cousin (\cite{SET2}).

For further results, the reader is referred to
\cite{Ba}, \cite{Koi} and \cite{RSE2} for the Euclidean space and
\cite{NSE} and \cite{ST1} for the hyperbolic space.

Now, we describe our results for $H$-surfaces in
$\hd \times \R$.

First, assume $M$ is compact, immersed into $\hd \times \R$,
with mean
curvature $1/2$, with boundary a round circle in a horizontal slice. Then $M$ is
part of the
rotational surface with vanishing Abresch-Rosenberg holomorphic quadratic
differential $Q$. An analogous result holds if $M$ has constant mean
curvature less than $1/2$ (Theorem \ref{oneboundary} and
\ref{T.Q=0}). In both situations $M$ is part of an entire rotational
vertical graph. More
generally, we prove that if $\partial M$ is a horizontal curve with curvature
greater than
one, then $M$ is a vertical graph, in particular $M$ has genus zero and inherits
the symmetries of its boundary (Theorem \ref{oneboundary}).  Assume now $M$ is
compact, embedded and has constant mean curvature $\leq 1/2$. Assume also that
$\partial M=C_1 \cup C_2$, where $C_1$ and $C_2$ are parallel horizontal curves
with curvature greater than one and $\dist (C_1,C_2)\geq \pi$.
Then $M$ inherits the symmetries of $C_1 \cup C_2$. Consequently, if $C_1$ and
$C_2$ are two parallel circles then $M$ is part of an embedded complete
rotational $H$-annulus. The last assertion follows from the geometric
classification of rotational $H$-surfaces (Section \ref{Sec.Rot}).

Finally, we consider surfaces $M$ in $\hd \times \R$ which are
regular up
to the asymptotic boundary and whose mean curvature is bounded away from zero.
In
\cite{Sa} R. Sa Earp gives examples suggesting that the
asymptotic boundary of such surfaces must lie on a
vertical line in $\pain (\hd) \times \R$. We prove that if the asymptotic
boundary
$\pain M$ in $\pain (\hd) \times \R$ is $C^1$
and the surface is $C^1$ up to $\pain M$, then each component of $\pain M$
is part of a vertical line. As far as we know, the same question is open, if
one relax the regularity up to the asymptotic boundary.
Observe that in the minimal case many
other possibilities can
occur (\cite{H}, \cite{Ne-Ro1}, \cite{Sa} and \cite{ST}). In
particular, in
\cite{Ne-Ro1} B. Nelli and H. Rosenberg solve the Dirichlet Problem in
$\hd \times \R$ for any
Jordan curve $\gamma$ in the asymptotic boundary $\pain (\hd) \times \R$
that is a vertical graph.

The paper is organized as follows. Section \ref{S.2} and \ref{S.3} deal with
compact $H$-surfaces with boundary either one curve or two curves, with
curvature
greater than one. In Section \ref{S.line}, we study the behavior
of the asymptotic boundary of a $H$-surface with strictly positive mean
curvature. Finally, in the Appendix we discuss the geometry of rotational
$H$-surfaces for any $H\in \R$, since we will need to use these surfaces
throughout
the whole paper. Rotational surfaces in
$\hd\times \R$ have been studied in \cite{Ab-Ro}, \cite{H-H},
\cite{Mo-On}, \cite{Pe-Ri} and \cite{ST}.
Further results on $H$-surfaces in $\hd \times \R$ are in
\cite{HST}, \cite{Ne-Ro2}.

\section{$H$-surfaces with boundary a curve with curvature greater than one}
\label{S.2}

In this section we  discuss existence and uniqueness of compact \newline
$H$-surfaces, $|H|\leqs 1/2,$ with
boundary a curve on a horizontal slice, with curvature greater than
one. Then, we study under which conditions a compact $H$-surface with
boundary a planar curve inherits the
symmetries of its boundary. In particular we  deal with the circle
boundary case.

Let $H\in (0, 1/2]$. Denote by $t$ the third coordinate in $\hd \times \R$ and
by $\sigma$ the origin in $\hd$(in the disk model of
$\hd $ we have $\sigma =0$).
We denote by $S^H$ the simply connected embedded
surface in
$\hd \times \R$ with constant mean curvature $H$, invariant by rotation about
the axis $\{\sigma \} \times ~\R$, tangent to
$\hd \times \{0\}$ at $\sigma$.  Recall that $S^H$ is an
entire vertical graph (Section \ref{Sec.Rot}).

Let us recall the Convex Hull Lemma (\cite{BR}).  Let $K$ be a
compact set in $\h^2\times\re.$ For any $H\in(0,\frac{1}{2}],$ we
define ${\mathcal F}^H_K $ as follows. A surface $B$ belongs
to ${\mathcal F}^H_K $ if $K$ is contained  in the mean
convex side of $B$ and if it is
obtained from $S^H$
either by vertical and horizontal translations or by  symmetry with
respect to a horizontal slice.

Using the maximum principle one can prove the following.

\begin{lemma} {\bf (Convex Hull Lemma)~\cite{BR}}
\label{convexhull}
 Let  $M$ be  a compact surface immersed in  $\h^2\times\re$ with constant mean curvature $H\in(0,\frac{1}{2}].$
Then $M$ is contained in the convex hull of the family ${\mathcal
F}^H_{\partial M}.$
\end{lemma}

We recall the Flux Formula (\cite{BSE1}, Appendix B in
\cite{BS} and \cite{HLR}). Let $\Sigma$  and $Q$ be two  compact,
smooth, not necessarily connected surfaces in $\h^2\times\re$ such
that their  boundaries coincide. Assume that there exists a compact
domain $U$  in $\h^2\times\re$  such that the boundary of $U$ is
$\partial U=\Sigma\cup Q$ and it is orientable. Notice that the
boundary of $U$ is smooth except perhaps on $\partial
\Sigma=\partial Q.$

 Let  $n_{\Sigma},$ $n_Q$ be the unit normal fields to
$\Sigma$ and $Q$ respectively that point inside   $U.$ Denote by $\nu$ the unit
conormal to $\Sigma$ along $\partial \Sigma,$ pointing inside
$\Sigma.$ Finally assume  that $\Sigma$ is  a compact surface with
constant mean curvature $H.$ Let $Y$ be a Killing vector field in
$\h^2\times\re.$ One can prove the {\it Flux Formula}

\begin{equation}\label{fluxformula}
 \int_{\partial\Sigma}\!\!\!
\langle Y,\nu\rangle=2H\!\!\int_Q\!\! \langle Y,n_Q\rangle.
\end{equation}

In $\h^2\times\re$ there is a natural  notion of  {\em vertical
graph} ( ~\cite{HST},  ~\cite{BR} or ~\cite{Sa}): let $\Omega$ be a subset
of
$\h^2$ and let $u:\Omega\lra\re$ be a $C^2$ function. The
{\it vertical graph} of $u$ is the subset of
 $\h^2\times\re$ given by
\begin{equation}
\nonumber \left\{(x,y,t)\in\Omega\times\re \ \vert \
t=u(x,y)\right\}.
\end{equation}

We choose the unit normal vector field to the graph of $u$
with positive third component  and we compute  the mean curvature
(of the graph) with respect to  it.

The graph of a function $u:\h^2\lra\re$ has the function $H$ as mean
curvature if and only if  $u$ satisfies the following partial differential
equation.

\begin{equation}\label{graph}
\mathrm{div}_{\h}\!\! \left(\frac{\nabla_{\h}
u}{W_u} \right)\!=\! 2H,
\end{equation}
where $\mathrm{div}_{\h},$ $\nabla_{\h} $ are the hyperbolic divergence and
gradient respectively and $W_u=\sqrt{1+|\nabla_{\h}u|_{\h}^2},$
being  $|\cdot|_{\h}$ the norm in  $\h^2.$

Consider the halfspace model for $\h^2,$ with Euclidean
coordinates $x,y,$  $y>0.$ In this model, equation (\ref{graph})
takes the following form

\begin{equation}\label{graphhalfspace}
\mathrm{div}\!\!\left(\frac{\nabla u}{W_u}\right)\!=\!\frac{2H}{y^2},
\end{equation}
where $\mathrm{div}$ is the Euclidean divergence  and
$W_u=\sqrt{1+y^2(u_x^2+u_y^2)}.$

The following Theorem is a consequence of a result proved by L. Hauswirth, H.
Rosenberg and
J. Spruck (\cite{HRS}, Th. 3.2). We sketch a proof for completeness.

\begin{theorem}\label{graphexistence}
 Let $\Omega$  be a  domain in
$\h^2$ such that $\partial\Omega$ is  $C^{2,\alpha}$ and has
curvature greater than one. For any $H\in [- \frac{1}{2},\frac{1}{2}]$, there
exists a vertical graph $G_H$ with constant mean curvature $H$ and
boundary $\partial\Omega.$
\end{theorem}

{\bf Proof.}

By formula (\ref{graphhalfspace}), one has to solve the following
Dirichlet problem

\begin{equation*}
\left\{\begin{array}{l}
\!\! \!F[u]=\mathrm{div}\left(\frac{\nabla u}{W_u}\right) -\frac{2H}{y^2}=0 \ \
{\rm
in} \ \  \Omega ,\\ \\
\!\! \! u=0  \ {\rm on } \ \partial\Omega .\\ \\
\end{array}\right.
\end{equation*}

In order to prove existence, we use the continuity method. For every
$t\in[0,1],$ consider the Dirichlet problem

\begin{equation}\label{dirichletfamily}
\left\{\begin{array}{l}
\!\!\! F^{t}[u]=\mathrm{div}\left(\frac{\nabla
u}{W_{u}}\right)-\frac{2Ht}{y^2}=0 \ \
{\rm in} \ \ \Omega , \\ \\
\!\!\! u=0  \ {\rm on } \ \partial\Omega .\\ \\
\end{array}\right.
\end{equation}

Let $S=\{t\in[0,1] \ \ |\ \ {\rm  there \ exists\ a \ solution \ of\
(\ref{dirichletfamily})}\}$. Observe that
$u\equiv 0$ is a solution of (\ref{dirichletfamily}) for $t=0$,
hence $0\in S$. If one proves that $S$ is open and closed, then
$1\in S$ and the desired solution  is a solution of
(\ref{dirichletfamily}) for $t=1.$

That $S$ is open follows from the Implicit Function Theorem.
$S$ closed follows from $C^{2,\alpha}$  {\em a-priori}
estimates for solutions of (\ref{dirichletfamily}). By Schauder's
theory (\cite{GT}), $C^{2,\alpha}$ {\em a-priori} estimates
follow from $C^1$ {\em a-priori} estimates. The Convex Hull Lemma
guarantees $C^0$  estimates and boundary gradient bounds on
solutions of (\ref{dirichletfamily}).
Therefore, we infer with  Theorem 3.1 in \cite{Sp}
that boundary gradient bounds
implies interior
gradient bounds. \qed

\begin{theorem}
\label{oneboundary}
 Let $M$ be a compact surface immersed  in
$\h^2\times \re$ with boundary a $C^{2,\alpha}$ Jordan curve  $C$ with
curvature greater than one, contained in the slice
$\hd \times \{0\}$. Assume
that $M$ has constant mean curvature $H\in(0,\frac{1}{2}].$ Then $M$
is a vertical graph (given by Theorem \ref{graphexistence}). In particular $M$
has genus zero and inherits
the same symmetries of its boundary. If $C$ is a circle, then $M$ is
a part of the simply connected rotational surface containing
$C$ with constant mean curvature $H$.
\end{theorem}

{\bf Proof.}

Denote by $P=\{t=0\}$ the slice containing the boundary curve $C$
and denote by $\Omega$ the domain in $P,$ bounded by $C.$ Consider the
convex hull of the family ${\mathcal F}^H_C.$ As $C$ has curvature
greater than one, for any point of $p\in C$ there is a surface of
the family ${\mathcal F}^H_C$  tangent  to $C$ at $p.$ By Lemma
\ref{convexhull}, $M$  is contained in the convex hull of the family
${\mathcal F}^H_C$, hence $M\backslash \partial M$ is entirely contained in
the vertical cylinder over $C$ and $M$ does not meet $P$ outside $\Omega$.

By Theorem \ref{graphexistence} there exists a graph $G_H$ on $\Omega,$
with boundary $C$ and constant mean curvature $H.$ We can choose
$G_H$ to be contained in the halfspace $t\geq 0$  with   mean
curvature vector pointing  downwards.
 We will  prove that $M$ is contained in one of the two halfspaces
determined by the slice $P.$

First assume that $M$ is embedded. Lift up  $G_H$ to be above  $M,$
then move   $G_H$ down   towards $M$: by the maximum principle,  one
can not touch $M$ till the boundary of $M$ and the boundary of $G_H$
coincide. Hence  $M$ lies below  $G_H.$ Then $G_H\cup M$ bounds a
domain $U$ in $\h^2\times \re$ and the mean curvature vector of $M$
points either inside $U$ or outside $U.$

Assume, by contradiction, that $M$ has points in both halfspaces \newline
$\{t> 0\},$ $\{t< 0\}.$ By the maximum principle, the mean
curvature vector at an highest point of $M$ points downwards, i.e.
outside $U,$  while the  mean curvature vector at a lowest point of
$M$ points upwards, i.e. inside $U.$ This is a contradiction. Then
$M$ is contained in a halfspace, say $\{t> 0\}$, $M\cup \Omega$
bounds a domain $W$ in  $\h^2\times \re$ and the mean curvature
vector of $M$ points  inside $W.$  We already know that $M$ lies
below  $G_H.$ Now move down $G_H$ to be disjoint from $M,$ then lift
$G_H$ up towards $M.$ At a first interior   interior contact point
the mean curvature vectors of $M$ and $G_H$  coincide. Contradiction
by the maximum principle. Then one can lift $G_H$ till the
boundaries of $M$ and $G_H$ coincide, that is $M$ is above $G_H.$
This implies that $G_H\equiv M,$ as desired.

Now, by applying  Alexandrov reflection method with vertical geodesic planes,
one obtains that $M$ has all the symmetries of $C.$ In particular,
if $C$ is a circle, then  $M$ is part of the simply connected
rotational surface $S^H$ containing $C$. Observe that if $M$ is embedded a
simpler alternative argument based on Alexandrov Reflection Principle, using the
horizontal slices, yields that $M$ is a vertical graph.

Now let  $M$ be   an immersed surface
(not necessarily embedded).

First assume  that $C$ is a circle.

By Lemma \ref{convexhull}, $M$ is contained in the convex hull of
the family ${\mathcal F}^H_C.$ As $C$ is a circle, this convex hull
is the domain bounded by the compact part of $S^H$ containing $C$,
say $B^1_C,$  and its symmetry with respect to the slice $P,$ say
$B^2_C.$

Let $\nu^1_3,$ $\nu_3,$ $\nu^2_3$ the third components of the
inward unit conormal along $C$ of $B^1_C,$ $M,$ $B^2_C$
respectively. As $M$ is between $B^1_C$ and $B^2_C,$ then at any
point of $C$

\begin{equation}\label{ineq}
 \nu^2_3\leq \nu_3\leq \nu^1_3.
\end{equation}

Consider the Flux Formula for $M$, with $\Omega$ equals to the planar
domain bounded by $C,$ $Y=(0,0,1)$  and $n_\Omega=(0,0,\pm 1)$ according
to the orientation given by $M.$  In order to  fix ideas, assume
that $n_\Omega=(0,0,1).$ Formula (\ref{fluxformula}) yields

\begin{equation}\label{eq}
 \int_C\!\!\nu_3= 2H Area(\Omega).
\end{equation}

Now, consider  the Flux Formula  for $B^1_C$ and $B^2_C,$ with $\Omega$
equal to the planar domain bounded by $C,$  $Y=(0,0,1).$  By
formula (\ref{fluxformula})

\begin{equation}\label{eq1}
\int_C\!\! \nu^1_3=2H Area(\Omega)=-\!\!\int_C\!\!\nu^2_3.
\end{equation}

Then,  equalities (\ref{eq}) and  (\ref{eq1}) yield

\begin{equation}
\label{eq2} -\!\int_C\!\! \nu^2_3\!=\!\int_C\!\! \nu_3 \!=\!\int_C\!\!\nu^1_3.
\end{equation}

If, in equation (\ref{ineq}), the inequalities are strict at  every
point of $C$ then, one has a contradiction by equation (\ref{eq2}).

Then,  there is at least one point $p$ in $C$ such that $\nu_3$
agrees with either $\nu^2_3$ or $\nu^1_3$ at $p$. Therefore, by the
boundary maximum principle, $M$ coincides with either $B^2_C$ or
$B^1_C.$

If the boundary of $M$ is an embedded curve $C$ with curvature greater
than one the proof is
analogous. It is enough to  replace the caps $B^i_C,$ $i=1,2$  by
the graph  $G_H$ of  mean curvature $H$  and boundary $C $
 and its symmetry with respect to the slice $P$. \qed

\begin{remark}
Notice that one can prove the analogous results for a surface $M,$
whose mean curvature function $H(x,y)$ satisfies for any
$(x,y)\in \ov{\Omega}$:
$0< \vert H(x,y) \vert \leq ~\frac{1}{2}$.
\end{remark}


\section{$H$-surfaces  with boundary two parallel
curves  with curvature greater than one}\label{S.3}

 We say that
$C_1$ and $C_2$ are {\em parallel curves}, if they are congruent up
to a vertical translation.

\begin{theorem}
\label{smalldistancetheorem}
 Let $M$ be a compact embedded surface in
$\hip^2\times \R,$ with boundary two parallel, embedded $C^{2,\alpha}$
curves
$C_a\subset P_a=\{ t=a\}$, and
$C_{-a}\subset P_{-a}=\{t=-a\},$  with curvature greater than
one. Assume that $M$ has constant mean curvature $H$, with $|H|\leqs
1/2.$ Then $M$ is symmetric with respect to the horizontal slice
$\{t=0\}$. If $2a \geq \pi,$
then $M$ is contained in the closed slab $\{-a\leqs t\leqs a\},$
with $M\cap (P_a\cup P_{-a})=C_a\cup C_{-a}.$ Furthermore, $M$
inherits the symmetries of $C_a\cup C_{-a}.$
\end{theorem}

{\bf Proof.}

Let $D_a$ be the bounded  domain in $P_a$ with boundary $C_a.$ Let
$\ext(D_a)=P_a\setminus D_a.$  By the Convex Hull Lemma,
$M\cap\ext(D_a)=C_a$ and $M\cap\ext(D_{-a})=C_{-a}.$ For any point
$p\in C_a,$ there exists a circle containing $D_a$ in its interior,
tangent to $C_a$ at $p.$ Let $Z_p$ be the cylinder over such a circle.
Denote by $Z_p^+$ the mean convex open domain of $\h^2\times \re$
bounded by $Z_p.$ Each $Z_p$ has mean curvature greater than
$\frac{1}{2}$ hence, by the maximum principle, $M$ is contained in
$\cap_{p\in C_a} Z_p^+.$ Notice that $\partial (\cap_{p\in C_a}
Z_p^+)$ is the cylinder over $C_a.$ We call it $Z$ and we have
that $M\cap Z=C_a\cap C_{-a}.$ Then, we apply Alexandrov reflection
with horizontal slices to infer that the slice $\{t=0\}$ is a plane of
symmetry for $M.$

Now, assume that $2a \geq \pi.$

By Theorem \ref{graphexistence}, there exists a graph $S_1$ with
boundary $C_a$  with mean curvature $H$ and mean curvature pointing
downward,  and a graph $S_2$ with boundary $C_{-a},$  with mean
curvature $H$ and mean curvature vector pointing upward. By the
maximum principle $M$ lies below $S_1$ and above $S_2,$ hence
 $M\cup S_1\cup S_2$
is a closed embedded surface, not smooth along $C_a\cup C_{-a}.$ Let
$U$ be the domain in $\h^2\times \re$ bounded by $M\cup S_1\cup
S_2.$ We claim that  the mean curvature vector at any point of  $M$
points towards $U.$  Consider the family of minimal catenoids
(Section \ref{Sec.Rot}). Let us recall their shape. For any $t\in (0,
\frac{\pi}{2})$ there exists a catenoid bounded by two circles at
infinity at height $\pm t.$ When $t\to 0,$ the catenoids
tend to the double covering of the slice $\{t=0\}$ and when $t
\to \frac{\pi}{2},$ the catenoids tend to infinity.
Then, one comes with catenoids from infinity towards $M.$ Let
$p$  be  the first contact point between $M$ and one of the
catenoids. By the maximum principle, the mean curvature vector of
$M$  at $p$ points  inside $U,$ then it points inside $U$ at any
point of $M.$

Now we prove that $M$ is contained in the closed slab
$\{-a\leqs t\leqs a\}$, with $M\cap \{t=\pm a\}=C_a\cup
C_{-a}.$

Assume by contradiction that $M$ has some points above the slice
$t=a,$ and let $q$ be a highest point of $M.$ As $S_1$ is above
$M,$ the mean curvature vector at $q$ points upwards. This gives a
contradiction by comparing $M$ with a horizontal slice. Now, assume
by contradiction that $M$ has some points below  the slice $t=-a,$
and let $q$ be a lowest point of $M.$ As $S_2$ is below $M,$ the
mean curvature vector at $q$ points downwards. And this gives a
contradiction by comparing $M$ with a horizontal slice.

By applying the Alexandrov reflection method with vertical geodesic planes, we
obtain that $M$ inherits the symmetries of $C_a\cup C_{-a}$. \qed

\begin{corollary}
\label{smalldistancecorollary} Let $M$ be a compact embedded surface
in $\h^2\times \re$ with boundary two parallel circles. Assume that
$M$ has constant mean curvature $H$, with $|H|\leqs 1/2.$  Let $d$ be
the distance of the two boundary curves. If $d\geq \pi,$ then $M$ is
part of an embedded complete rotational surface of constant mean
curvature $H.$
\end{corollary}

{\bf Proof.}

By the previous Theorem,  $M$ is part of a rotational surface  with
mean curvature $H.$  From the proof of Lemma \ref{smalldistancetheorem}, we
infer that the mean curvature vector of $M$ points  towards the
interior of $M\cup D_a\cup D_{-a},$ where $D_a\subset P_a,$
$D_{-a}\subset P_{-a}$ are the domains bounded by  $C_a$ and  $C_{-a}$
respectively. By the geometric classification of the rotational
surface with constant mean curvature $H,$ $M$ must be a part of an
embedded complete rotational surface. \qed

\section{Surfaces with nonempty asymptotic boundary}\label{S.line}

In \cite{Sa} R. Sa Earp describes many examples of complete
$H$-surfaces with $H>0$. When
the asymptotic boundary is nonempty, then it consists of parts of straight
lines.

\begin{theorem}
Let $M$ be a surface in $\hd\times \R$ with mean curvature satisfying
$0<\delta\leq H(p)$ at any point $p\in M$.
Assume that the asymptotic boundary of $M$ in $\pain (\hd) \times \R$ is a $C^1$
curve and that $M$ is
$C^1$ up to the asymptotic boundary (eventually $M$ has nonempty finite
boundary).

Then, each connected component of the asymptotic boundary of $M$ is part of a
vertical straight line in $\pain (\hd) \times \R$.
\end{theorem}

{\bf Proof.}   \par
Let $C$ be a connected component of the asymptotic boundary of $M$,
$C\subset \pain M\subset \pain (\hd) \times ~\R$.
We show that $C$ is vertical at the point $p$ for any
$p \in C$. We can assume that
$p=1\in \pain(\hd) \times \{0\}$.

Suppose that $C$ is not vertical at $p$. Then, an open neighborhood of $p$ in
$C$ is a graph $(e^{i\theta},t(\theta))$ for $\theta$ in an open interval
around $0$ ($p=1= e^{i 0}$) with $t(0)=0$. Let $\varepsilon >0$ be a small
number.
There exists a real number
$\nu=\nu (\varepsilon)>0$ such that $\vert t(\theta) \vert <\varepsilon/2$ for
any
$\theta \in [-\nu , \nu ]$. We set $p_1:=(e^{i\nu},t(\nu))$ and
$p_2:=(e^{-i\nu},t(-\nu))$ and we call $C_\var \subset C$ the closed subarc of
$C$ bounded by $p_1$ and $p_2$ and containing $p$
$$
C_\var\! \subset \! \pain (\hd) \!\times \! (-\var/2,\var/2),\
\partial C_\var =\{p_1,p_2\}.
$$

Since $M$ is $C^1$ up to the asymptotic boundary, there exists a simple arc
$C^\prime_\var \subset M$ with asymptotic boundary $p_1$ and $p_2$ such that
the connected component $M_\var$ of $M\backslash C^\prime_\var$ containing
$C_\var$ in its asymptotic boundary satisfies
$M_\var \subset \hd \times ( -\var,\var)$. Summarizing, we have
$$
M_\var \!\subset \! M\cap \hd \times ( -\var,\var ),\ \partial
M_\var=C^\prime_\var,\ \pain M_\var=C_\var \! \subset \! \pain (\hd) \times
(-\var/2,\var/2).
$$

Let $\Pi :\ov{\hd} \times \R \rightarrow \ov{\hd}$ be
the projection on the first two coordinates.

 For $\var$ small enough,
$\Pi(C_\var)\subset \pain \hd$ is an open arc with end points $\Pi(p_1)$ and
$\Pi(p_2)$ and containing $1=p$ in its interior. Moreover
$\Pi (M_\var) $
is an open subset of $\hd$ with asymptotic boundary $\Pi (C_\var)$.

Let $\Gamma\in \hd$ be the geodesic $(-1,1)$ and let $T_s, s>0$,
be the hyperbolic translation along $\Gamma$ defined as follows:
any point $x\in \Gamma$ is sent to the point of $\Gamma$ between $x$ and $-1$
whose hyperbolic distance from $x$ is $s$. Then, extend $T_s$ to
$\hd \times \R$ by vertical translation.

Let $\eta >0$ be a small number. For $s$ great enough and $\var$ small enough,
the curve $T_s(\Pi(C_\var^\prime))$ is inside the open euclidean ball
centered at $-1$ with radius $\eta$. Thus $T_s(\Pi(M_\var))$ is the
connected component of $\hd\backslash T_s(\Pi(C_\var^\prime))$ containing 1 in
the asymptotic boundary.

The vertical projection of the surface $T_s(M_\var)$ covers a large part
of $\hd$ in the euclidean sense.
 Its boundary is the simple arc
$T_s(C_\var^\prime)$. Let
$\gamma_s\subset \pain (\hd) \times  (-\var,\var)$ be a $C^1$
arc with end points $T_s(p_1)$ and $T_s(p_2)$ such that
$\gamma_s \cup T_s(C_\var)$ is a Jordan curve which projects one to one onto
$\pain (\hd) \times \{0\}$.

Finally, let $R_s \subset \hd \times \R$ be any embedded smooth disk, disjoint
from
the interior of $T_s(M_\var)$, with finite boundary
$T_s(C_\var^\prime)$ and asymptotic boundary $\gamma_s$:
\begin{align*}
 R_s \!\cap \! T_s(M_\var) &= T_s(C_\var^\prime),  \\
 \partial R_s &= T_s(C_\var^\prime),  \\
 \pain R_s &= \gamma_s .
\end{align*}
Then $R_s \cup T_s(M_\var)$ is an embedded simply connected surface, it
is $C^0$ along $T_s(C_\var^\prime)$ and smooth everywhere else. The surface
$R_s \cup T_s(M_\var)$ separates $\hd \times \R$ in two connected
components and its asymptotic boundary is the Jordan curve
$\gamma_s \cup T_s(C_\var) \subset \pain (\hd) \times [-\var,\var]$.

First, we assume that the mean curvature vector of
$ T_s(M_\var)$ points towards the connected component
containing $\hd \times (\var , +\infty)$.

We can assume that $\delta<1/2$.
Let us consider the simply connected $H$-surface $S^\delta$
given in Proposition
\ref{P.H<=1/2} with $H=\delta$ and with rotational axis equals to the vertical
geodesic $\{0\} \times \R$. Now, lift up $S^\delta$ to be above
$R_s \cup T_s(M_\var)$, then move $S^\delta$ down. By our construction and the
geometry of $S^\delta$ the first contact with $R_s \cup T_s(M_\var)$ will be at
an interior point of $T_s(M_\var)$. This gives a contradiction with the maximum
principle. If the mean curvature vector of $T_s(M_\var)$ points towards the
other component, one does the same reasoning with the surface obtained from
$S^\delta$ by symmetry with respect to the slice $\hd\times \{0\}$.
Therefore, the asymptotic boundary
$C$ is vertical at any point $p\in C$. \qed

\section{Appendix: Geometric behavior of the rotational $H$-surfaces in
$\hi2 \times \R$}\label{Sec.Rot}

In the Appendix we describe in details
the geometric behavior of rotational $H$-surfaces. Our discussion
is based on formulae founded in \cite{ST}. We recall that rotational
surfaces in
$\hd\times \R$ have been studied in \cite{Ab-Ro}, \cite{H-H},
\cite{Mo-On}, \cite{Pe-Ri} and \cite{ST}.

 We work with the disk model for $\hd $, so that
$$
\hd =\{(x,y)\in \m 2,\ x^2+y^2<1 \},
$$
and the metric is
$$
\rmd s_\mathbb{H} ^2=\left( \dfrac{2}{1-(x^2+y^2)}\right)^2
(\rmd x^2 + \rmd y^2).
$$
Therefore the product metric on $\hd \times \R$ reads as follows
$$
\rmd \tilde s ^2=\left( \dfrac{2}{1-(x^2+y^2)}\right)^2
(\rmd x^2 + \rmd y^2) + \rmd t^2,
$$
where $(x,y)\in \hd$ and $t\in \R$.
We consider the following particular geodesic of $\hd$
\begin{equation*}
\Gamma \!= \!\{(x,0), \ x\in (-1,1)\,\}\!\subset \! \hd .
\end{equation*}

Up to ambient isometry,
we can assume the rotational surfaces
are generated by curves in
the vertical geodesic plane $P= \Gamma \times \R \subset \hd\times\R$ and
that
the
rotational axis is the vertical geodesic $R:=\{(0,0)\}\times \R$.

On the geodesic $\Gamma$ we denote by $\rho \in \R$ the signed
distance to the origin $(0,0)$, thus $x=\tanh \rho/2$.
Therefore the metric on $P$ is
$$
\rmd s^2=\left( \dfrac{2}{1-x^2}\right)^2\rmd x^2 +\rmd t^2
=\rmd\rho^2 +\rmd t^2.
$$
Let us consider a curve in $P$ which is a vertical graph:
$c(\rho)=(\rho, \lambda (\rho))$ where $\lambda$ is a smooth real function
defined for $\rho \geq 0$. On the rotational surface generated by $c$ we
consider the orientation given by the unit normal field pointing up.
It is shown in \cite{ST} (formula (21)) that the curve $c$ generates a
rotational surface with constant mean curvature $H$ if and only if the function
$\lambda$ is given by
\begin{equation}\label{Eq.Graph.1}
\lambda(\rho)= \int ^{\rho}_\ast \frac{d+2H\ch r}{\sqrt{\sh^2 r - ( d+2H\ch
r)^2}}dr,
\end{equation}
where $d$ is a real parameter and $\ast$ is the minimum such that the condition
$\sh^2 r - ( d+2H\ch r)^2 \geq 0$ is satisfied. Using the isometry
$(x,y,t) \mapsto (x,y,-t)$ we can assume that $H\geq 0$.
We will analyze consecutively the cases $H=0, H\in (0,1/2]$ and $H>1/2$.

\bigskip

\begin{proposition}(Minimal rotational surfaces)\label{P.Min}
For each $d\geq 0$ there exists a complete minimal rotational
surface $\mathcal{M}_d$. The surface $\mathcal{M}_0$ is the
horizontal slice $\{t=0\}$. For $d>0$ the rotational surface
$\mathcal{M}_d$ (called catenoid) is embedded and homeomorphic to an
annulus. The distance between the rotational axis and the ``neck''
of $\mathcal{M}_d$ is $\argsh d$. The asymptotic boundary of
$\mathcal{M}_d$ is two horizontal circles in $\pain (\hd) \times \R$
and the vertical distance between them is a nondecreasing function
$h(d)$ satisfying $\lim_{d \to 0}h(d)=0$ and $\lim_{d\to
+\infty}h(d)=\pi$. Therefore $\mathcal{M}_d$ converges to the double
covering of the slice $\{t=0\}$ when $d$ goes to 0.

Moreover any minimal rotational surface is, up to an ambient iso\-me\-try, a
part of a complete surface $\mathcal{M}_d$.
\end{proposition}

\bigskip
{\bf Proof.}   \par
If the graph of a function $\lambda$ generates a minimal surface
we deduce from formula (\ref{Eq.Graph.1}) that
\begin{equation*}
\lambda (\rho)=  \int_{\argsh (d)} ^{\rho} \frac{d}{\sqrt{\sh^2 r - d^2}}dr.
\end{equation*}
Thus we have $\lambda \equiv 0$ for $d=0$.
For $d>0$ the function $\lambda $ is defined for
$\rho \geq \argsh (d) >0$ and the graph $c$ has a vertical tangent at
$\rho =\argsh d$. It is clear that $\lim_{\rho \to +\infty}\lambda(\rho)$
exists and is finite. Let us call $\wt{c}$ the union of
the curve $c$ with its
symmetry with respect to the horizontal geodesic $\{y=0\}$ of $P$.
Therefore  $\wt{c}$ is a complete curve which generates a rotational complete
and embedded minimal surface homeomorphic to an annulus.

As in formula (41) of \cite{ST}  we introduce the new coordinate $s$
setting $ds=\sqrt{1+\lambda^{\prime 2}} d\rho$ ($s$ is the arclength of the
graph $c(\rho)=(\rho,\lambda (\rho))$). Then using the formulae (49),
(36) and (37) of \cite{ST} we get
\begin{align*}
 \rho (s) &= \int_0^s \frac{\sh t}{\sqrt{(1+d^2)\ch^2 t -1}} dt +\argch
\sqrt{1+d^2}, \\
  &= \argch (\sqrt{1+d^2} \ch s).
\end{align*}
 and
\begin{equation*}
 \lambda \!\circ \!\rho (s)\!=\!\int_0^s\!\! \frac{d}{\sqrt{(1+d^2)\ch^2 t -1}}
dt,
\end{equation*}
for $s \geq 0$. Therefore we have
\begin{equation*}
 h(d)\!=\!2 \!\int_0^{+\infty}\!\!\! \frac{d}{\sqrt{(1+d^2)\ch^2 t -1}} dt.
\end{equation*}
Consider the positive function defined by
$f(t,d)=\frac{d}{\sqrt{(1+d^2)\ch^2 t-1}}$
for $t, d \geq 0$. Let $d_0 >0$ be any positive real number. Clearly the
integral $h(d)$ is convergent for any $d\in [d_0,+\infty )$. Moreover we have
\begin{equation*}
 \frac{\partial}{\partial d}f(t,d)=\frac{\sh t}{((1+d^2)\ch^2 t -1)^{3/2}}\cdot
\end{equation*}
We deduce that the integral
\begin{equation*}
 \int_0^{+\infty} \frac{\partial}{\partial d}f(t,d) dt,
\end{equation*}
is uniformly convergent for $d\geq d_0 >0$. Consequently the function $h(d)$
is differentiable on $[d_0,+\infty )$ and
\begin{align}
 h^\prime (d) &=  2 \int_0^{+\infty} \frac{\partial}{\partial d}f(t,d) dt ,\\
  &= 2 \int_0^{+\infty} \frac{\sh t}{((1+d^2)\ch^2 t -1)^{3/2}}dt
.\label{hauteur}
\end{align}
As this is true for $d\geq d_0>0$ for any $d_0>0$, then
$h$ is differentiable for $d>0$ and its derivative is given by
(\ref{hauteur}). We deduce that $h(d)$ is nondecreasing and we have
$\lim_{d \to 0}h(d)=0$.

Finally, from the inequalities
\begin{equation*}
 \frac{d}{\sqrt{(1+d^2)}}\frac{1}{\ch t} \leq
\frac{d}{\sqrt{(1+d^2)\ch^2 t-1}}  \leq \frac{1}{\ch t} ,
\end{equation*}
for any $d,t >0$, we get
\begin{align*}
 \lim_{d \to +\infty}h(d) &= 2 \int_0^{+\infty} \lim_{d \to +\infty}
\frac{d}{\sqrt{(1+d^2)\ch^2 t-1}} dt ,\\
  &= 2 \int_0^{+\infty}  \frac{dt}{\ch t} ,\\
  &= 2 \int_0^{+\infty}  \frac{du}{u^2+1},\ \ (u=\sh t) ,\\
  &= \pi .
\end{align*}
This concludes the proof. \qed

\bigskip

For later use we define the functions $g(\rho)$ and $f(\rho)$ setting for
$d\in\R$ and $H>0$,
\begin{align*}
 g(\rho) &= d+2H\ch \rho ,\\
 f(\rho) &= \sh^2\rho -(d+2H\ch \rho)^2 ,\\
 &= (1-4H^2)\ch^2\rho -4dH \ch \rho -1-d^2 ,
\end{align*}
so that $\lambda^\prime (\rho)=g(\rho)/\sqrt{f(\rho)}$.

\bigskip

\begin{lemma}\label{H<1/2}
Assume  $0<H< 1/2$.
We have $f(\rho)\geq 0$ if and only if
$\ch \rho\geq \frac{2dH+\sqrt{1-4H^2+d^2}}{1-4H^2}$. Let
$\rho_1\geq 0$ such that $\ch\rho_1= \frac{2dH+\sqrt{1-4H^2+d^2}}{1-4H^2}$,
then $f(\rho_1)=0$ and
$\rho_1=0$ if and only if $d=-2H$.
\begin{enumerate}
\item If $d> -2H$, then $\frac{-d}{2H}< \ch \rho_1$. Consequently
the function $\lambda$ is nondecreasing for $\rho \geq \rho_1>0$ and has a
nonfinite derivative at $\rho_1$.

\item If $d=-2H$, then
$\lambda^\prime(\rho)=\frac{2H\sqrt{\ch\rho -1}}
{\sqrt{(1-4H^2)\ch\rho +4H^2+1}}$. Therefore
the function $\lambda$ is defined for $\rho \geq 0$, it has a zero
derivative at $0$ and is nondecreasing for $\rho>0$.

\item If $d < -2H$, then there exists $\rho_0 > \rho_1 >0$ such that
$\frac{-d}{2H}=\ch \rho_0$. Consequently the function $\lambda$ is defined for
$\rho \geq \rho_1>0$ with a nonfinite derivative at $\rho_1$, it is
nonincreasing for \newline
$\rho_1<\rho<\rho_0$, has a zero deri\-vative at
$\rho_0$ and
it is nondecreasing for $\rho>\rho_0$.

\item For any $d$ we have $\lim_{\rho\to
+\infty}\lambda(\rho)=+\infty$.
\end{enumerate}
\end{lemma}

Next Lemma, is analogous to Lemma \ref{H<1/2} in the case $H=1/2$.
We observe that, in this case, the
set $\{\rho > 0 \mid f(\rho)>0\}$ is
nonempty if and only if $d<0$.

\begin{lemma}\label{H=1/2}
Assume $H=1/2$ and $d<0$. Then
$f(\rho)\geq 0$ if and only if $\ch \rho \geq\frac{1+d^2}{-2d}$.
Let $\rho_1\geq 0$ such that $\ch\rho_1= \frac{1+d^2}{-2d}$,
then $f(\rho_1)=0$ and
$\rho_1=0$ if and only if $d=-1$.
\begin{enumerate}
\item If $d\in ( -1,0)$, then $\frac{-d}{2H}< \ch \rho_1$. Consequently
the function $\lambda$ is nondecreasing for $\rho \geq \rho_1>0$ and has a
nonfinite derivative at $\rho_1$.

\item If $d=-1$, then
 $\lambda^\prime(\rho)=\frac{1}{\sqrt{2}} \sqrt{\ch\rho -1}$.
Therefore
the function $\lambda$ is defined for $\rho \geq 0$, it has a zero
deri\-vative at $0$ and is nondecreasing for $\rho>0$.

\item If $d < -1$ there exists $\rho_0 > \rho_1 >0$ such that
$\frac{-d}{2H}=\ch \rho_0$. Consequently the function $\lambda$ is defined for
$\rho \geq \rho_1>0$ with a nonfinite derivative at $\rho_1$, it is
nonincreasing for $\rho_1<\rho<\rho_0$, has a zero deri\-vative at
$\rho_0$ and
it
is nondecreasing for $\rho>\rho_0$.

\item For any $d$ we have $\lim_{\rho\to +\infty}\lambda(\rho)=+\infty$.
\end{enumerate}
\end{lemma}

The proof of Lemma \ref{H<1/2} and \ref{H=1/2} is a straightforward
computation taking into account Formula (\ref{Eq.Graph.1}).
As a consequence of Lemma \ref{H<1/2} and \ref{H=1/2}
we have the following results.

\begin{proposition}(Rotational $H$-surfaces with $\vert H\vert \leq1/2$)
\label{P.H<=1/2}

Assume $0<H\leq 1/2$. There exists
a one-parameter family $\mathcal{H}_d$, $d\in \R$ for $H<1/2$ and $d<0$ for
$H=1/2$, of complete rotational
$H$-surfaces.
\begin{enumerate}
\item For $d> -2H$, the surface $\mathcal{H}_d$ is a properly embedded annulus, symmetric with respect to the slice $\{t=0\}$, the distance
between the ``neck'' and the rotational axis $R=\{(0,0)\} \times \R$
is $\argch(\frac{2dH+\sqrt{1-4H^2+d^2}}{1-4H^2})$ for $H<1/2$ and
$\argch (\frac{1+d^2}{-2d})$ for $H=1/2$.

\item For $d=-2H$, the surface $\mathcal{H}_{-2H}$ is an entire vertical graph,
denoted by $S^H$. Moreover $S^H$ is contained in the halfspace
$\{t\geq 0\}$  and it is tangent to the slice $\hd\times \{0\}$ at
the point $(0,0,0)$.

\item For $d < -2H$, the surface $\mathcal{H}_d$ is a properly immersed
(and nonembedded) annulus, it is symmetric with respect to the slice
$\{t=0\}$, the distance between the ``neck'' and the rotational axis
$R$ is $\argch (\frac{2dH+\sqrt{1-4H^2+d^2}}{1-4H^2})$ for $H<1/2$
and $\argch (\frac{1+d^2}{-2d})$ for $H=1/2$.

\item In each of the previous case the surface is unbounded in the
$t$-coordinate. When $d$
tends to $-2H$ with either $d> -2H$ or
$d < -2H$, then the surfaces
$\mathcal{H}_d$
tends towards the union
 of $S^H$ and its symmetry with respect
to the slice $\{t=0\}$.
Furthermore, any rotational $H$-surface with $0<H\leq 1/2$ is, up to an ambient
isometry, a part of a surface of the family $\mathcal{H}_d$.
\end{enumerate}
\end{proposition}

\bigskip
{\bf Proof.}   \par
The result is a straightforward consequence  of Lemma \ref{H<1/2} and
\ref{H=1/2}. For
$d = -2H$, $\mathcal{H}_{-2H}$ is the rotational surface generated by the graph
of the function  $\lambda$.

For $d\not= -2H$, let $\gamma$ be the union of the graph of $\lambda$ joint
with its symmetry with respect to the slice $\{t=0\}$. Then
$\mathcal{H}_d$ is the rotational surface generated by the curve $\gamma$. \qed



\bigskip

Observe that, for $H>1/2$, the set $\{\rho > 0 \mid f(\rho)>0\}$
is nonempty
if and only if $d<-\sqrt{4H^2-1}$.

\bigskip

\begin{lemma}\label{H>1/2}
Let $H$ and $d$  satisfying $H>1/2$ and
$d<-\sqrt{4H^2-1}$.
 Then,
there exist two numbers
$0\leq \rho_1<\rho_2$ such that $\ch\rho_1=\frac{2dH+\sqrt{1-4H^2+d^2}}{1-4H^2}$
and
$\ch\rho_2=\frac{2dH-\sqrt{1-4H^2+d^2}}{1-4H^2}$. Therefore,
$f(\rho)>0$ if and only if $\rho_1<\rho<\rho_2$ and
$f(\rho_1)=f(\rho_2)=~0$.
\begin{enumerate}
\item If $d<-2H$, then $\rho_1>0$ and there exists a unique number
$\rho_0 \in (\rho_1,\rho_2)$ satisfying $g(\rho_0)=0$. Furthermore $g\leq 0$
on $[\rho_1,\rho_0)$ and  $g\geq 0$ on $(\rho_0,\rho_2]$. Consequently, the
function $\lambda$ is defined on $[\rho_1,\rho_2]$, has a nonfinite
derivative at $\rho_1$ and $\rho_2$, has a zero derivative at $\rho_0$,
is nonincreasing on $(\rho_1,\rho_0)$ and nondecreasing on $(\rho_0,\rho_2)$.

\item If $d=-2H$, then $\rho_1=0$ and
$\lambda^\prime(\rho)=\frac{2H\sqrt{\ch\rho -1}}
{\sqrt{(1-4H^2)\ch\rho +4H^2+1}}$.
Consequently, the function $\lambda$ is defined on $[0,\rho_2]$, is
nondecreasing, has
a zero derivative at 0 and a nonfinite derivative at $\rho_2$.

\item If $-2H<d<-\sqrt{4H^2-1}$, then $\rho_1>0$ and $g\geq 0$
on $[\rho_1,\rho_2]$. Therefore the function $\lambda$ is defined on
$[\rho_1,\rho_2]$, is nondecreasing and has nonfinite derivative at $\rho_1$
and $\rho_2$.
\end{enumerate}
\end{lemma}

The proof of Lemma \ref{H>1/2} follows from formula (\ref{Eq.Graph.1}) by
a computation. In the following Proposition we assume the notations of
Lemma \ref{H>1/2}

\begin{proposition}(Rotational surfaces with $H>1/2$)\label{P.H>1/2}

Assume $H>1/2$.
There exists a one-parameter family $\mathcal{D}_d$
of complete rotational $H$-surfaces, $d\leq -\sqrt{4H^2-1}$.
\begin{enumerate}
\item For $d<-2H$, the surface $\mathcal{D}_d$
is an immersed (and nonembedded) annulus, invariant by a vertical
translation and is contained in the closed region bounded by the two
vertical cylinders $\rho=\rho_1$ and $\rho=\rho_2$. Furthermore
$\rho_1 \to +\infty$ and $\rho_2 \to +\infty$ when $d\to -\infty$
and $\rho_1 \to 0$ and $\rho_2 \to\argch(\frac{4H^2+1}{4H^2-1})$
when $d\to -2H$. Such surfaces are analogous to the nodoids of
Delaunay in $\m 3$.

\item For $d=-2H$, the surface $\mathcal{D}_{-2H}$ is an embedded sphere
and the maximal distance from the rotational axis is
$\rho_2=\argch(\frac{4H^2+1}{4H^2-1})$.

\item For $-2H<d<-\sqrt{4H^2-1}$, the surface $\mathcal{D}_d$ is
an embedded annulus, invariant by a vertical translation and is
contained in the closed region bounded by the two vertical cylinders
$\rho=\rho_1$ and $\rho=\rho_2$. Furthermore $\rho_1 \to 0$ and
$\rho_2 \to\argch(\frac{4H^2+1}{4H^2-1})$ when $d\to -2H$ and both
$\rho_1, \rho_2 \to \argch(\frac{2H}{\sqrt{4H^2-1}})$ when $d\to
-\sqrt{4H^2-1}$. Moreover
$\rho_1<\argch(\frac{2H}{\sqrt{4H^2-1}})<\rho_2$. Such surfaces are
analogous to the undoloids of Delaunay in $\m 3$.

\item For $d=-\sqrt{4H^2-1}$, the surface $\mathcal{D}_{ -\sqrt{4H^2-1}}$
is the vertical
cylinder over the circle with hyperbolic radius
$\argch(\frac{2H}{\sqrt{4H^2-1}})$.
\end{enumerate}
\end{proposition}

\bigskip

{\bf Proof.}   \par
For $d=-\sqrt{4H^2-1}$ we get the limit case of a vertical cylinder given by
$\ch\rho=\frac{2H}{\sqrt{4H^2-1}}$. A straightforward computation shows the
mean curvature of such a cylinder is $H$.

For the other cases the proof is a straightforward consequence of Lemma
\ref{H>1/2}. Let $\gamma$ be the union of the graph of $\lambda$ joint with its
symmetries with respect to the horizontal slices on which $\lambda$ is
vertical. When $d=-2H$, $\gamma$ is a compact arc and for $d\not=-2H$, $\gamma$
is periodic and complete, embedded only when $d>-2H$. \qed

\bigskip

From Proposition \ref{P.Min}, \ref{P.H<=1/2}, \ref{P.H>1/2}
and Proposition 26 in \cite{ST}, we infer the following
classification of
rotational $H$-surfaces with vanishing Abresch-Rosenberg holomorphic quadratic
differential. The classification of $H$-surfaces with vanishing
Abresch-Rosenberg holomorphic quadratic differential is established in
\cite{Ab-Ro}.

\begin{theorem}\label{T.Q=0}
Let $M$ be a rotational $H$-surface, $H\geq 0$, with vanishing
Abresch-Rosenberg holomorphic quadratic differential. We have up to congruence:
\begin{enumerate}
\item If $H=0$ then $M$ is a slice $\hd \times \{t\}$.
\item If $H>1/2$ then $M$ is an embedded two-sphere.
\item If $H=1/2$ then $M$ is the entire vertical graph $S^{1/2}$.
\item If $H<1/2$ then $M$ is either the entire vertical graph $S^H$ or the
embedded annulus $\mathcal{H}_{2H}$.
\end{enumerate}
\end{theorem}

\end{document}